\documentclass[12pt,psamsfonts,leqno,oneside,letterpaper]{amsart}

\usepackage{color}

\theoremstyle{definition}
\newtheorem{para}{}[section]
\newtheorem{subpara}{}[para]
\newtheorem{remark}[para]{Remark}
\newtheorem{remarks}[para]{Remarks}
\newtheorem{notation}[para]{Notation}
\newtheorem{convention}[para]{Convention}
\newtheorem{definition}[para]{Definition}
\newtheorem{definitions}[para]{Definitions}

\newcommand\Alternatives{\begin{enumerate}[(i)]}
\newcommand\EndAlternatives{\end{enumerate}}
\newcommand\Conditions{\begin{enumerate}[(1)]}
\newcommand\EndConditions{\end{enumerate}}

\theoremstyle{plain}
\newtheorem{theorem}[para]{Theorem}
\newtheorem{lemma}[para]{Lemma}
\newtheorem{proposition}[para]{Proposition}
\newtheorem{corollary}[para]{Corollary}
\newtheorem{conjecture}[para]{Conjecture}
\newtheorem{claim}[subpara]{}

\numberwithin{equation}{para}
\numberwithin{figure}{section}

\newcommand\Number{\begin{para}}
\newcommand\EndNumber{\end{para}}
\newcommand\Definition{\begin{definition}}
\newcommand\EndDefinition{\end{definition}}
\newcommand\Definitions{\begin{definitions}}
\newcommand\EndDefinitions{\end{definitions}}
\newcommand\Theorem{\begin{theorem}}
\newcommand\EndTheorem{\end{theorem}}
\newcommand\Conjecture{\begin{conjecture}}
\newcommand\EndConjecture{\end{conjecture}}
\newcommand\Remark{\begin{remark}}
\newcommand\EndRemark{\end{remark}}
\newcommand\Remarks{\begin{remarks}}
\newcommand\EndRemarks{\end{remarks}}
\newcommand\Convention{\begin{convention}}
\newcommand\EndConvention{\end{convention}}
\newcommand\Notation{\begin{notation}}
\newcommand\EndNotation{\end{notation}}
\newcommand\Lemma{\begin{lemma}}
\newcommand\EndLemma{\end{lemma}}
\newcommand\Proposition{\begin{proposition}}
\newcommand\EndProposition{\end{proposition}}
\newcommand\Corollary{\begin{corollary}}
\newcommand\EndCorollary{\end{corollary}}
\newcommand\Claim{\begin{claim}}
\newcommand\EndClaim{\end{claim}}
\newcommand\Proof{\begin{proof}}
\newcommand\EndProof{\end{proof}}
\newcommand\Equation{\begin{equation}}
\newcommand\EndEquation{\end{equation}}

\newcommand\Bullets{\begin{itemize}}
\newcommand\EndBullets{\end{itemize}}

\newcommand\tM{\widetilde M}

\newcommand\tGamma{\widetilde \Gamma}

\newcommand\ZZ{{\mathbb Z}}

\newcommand\NN{{\mathbb N}}

\newcommand\QQ{{\mathbb Q}}
\newcommand\HH{{\mathbb H}}

\newcommand\Hom{\mathop{\rm Hom}}

\newcommand\Fix{\mathop{{\rm Fix}}}

\newcommand\isomplus{\mathop{{\rm Isom}_+}}
\newcommand\isom{\mathop{{\rm Isom}}}

\begin{document}
\author{Peter B. Shalen}
\address{Department of Mathematics, Statistics, and Computer Science (M/C 249)\\  University of Illinois at Chicago\\
  851 S. Morgan St.\\
  Chicago, IL 60607-7045} \email{shalen@math.uic.edu}
\thanks{Partially supported by NSF grant DMS-0906155}
\dedicatory{Dedicated to Bus Jaco on the occasion of his 70th birthday}

\title{A Generic Margulis Number for Hyperbolic $3$-Manifolds}

\begin{abstract}
  We show that $0.29$ is a Margulis number for all but finitely many
hyperbolic 3-manifolds. The finitely many exceptions are all closed.
\end{abstract}

\maketitle

\section{Introduction}

If $M$ is an orientable hyperbolic $3$-manifold, we may write
$M=\HH^3/\Gamma$, where 
$\Gamma\cong\pi_1(M)$ is a  discrete,
torsion-free subgroup of $\isomplus(\HH^3)$, uniquely determined
up to conjugacy by the hyperbolic structure of $M$.

\begin{definition} Suppose that $M=\HH^3/\Gamma$ is a non-elementary orientable hyperbolic $3$-manifold. A {\it
Margulis number} for $M$ (or for $\Gamma$) is a number $\mu>0$ such that:
\begin{itemize}
\item
If $P\in \HH^n$, and if $x$ and $y$ are elements of $\Gamma$ such that $\max(d(P,x\cdot P),d(P,y\cdot P))<\mu$,   then $x$ and $y$ commute.
\end{itemize}
\end{definition}

Here, and throughout the paper, $d$ denotes  hyperbolic distance on $\HH^3$.

Note that if $\mu$ is a Margulis number for $M$ then any number in the interval $(0,\mu)$ is also a Margulis number for $M$. Note also that if $\Gamma$ is abelian then {\it any} positive number is a Margulis number for $M$.

The main result of this paper, which will be proved in Section \ref{proof section}, is:

\begin{theorem}\label{da goods}
Up to isometry there are at most finitely many orientable hyperbolic $3$-manifolds for which $0.29$ is not a Margulis number. Furthermore, all such manifolds are closed.
\end{theorem}

The first assertion of this theorem may be expressed as saying that $0.29$ is a ``generic Margulis number'' for orientable hyperbolic $3$-manifolds.

To put Theorem \ref{da goods} in context, I will first discuss a generalization of the definition of Margulis number that I gave above.  If $M$ ia a hyperbolic $n$-manifold, and if we write
$M=\HH^n/\Gamma$ where $\Gamma\le\isom(\HH^n)$ is discrete and
torsion-free, 
one may define a Margulis number for $M$ (or for $\Gamma$) to be a constant $\mu>0$ such that for every
$P\in\HH^n$, the elements $x\in\Gamma$ such
that $d(P,x\cdot P)<\mu$ generate a subgroup which is {\it elementary} in the sense that is
 has an abelian subgroup of finite index. If $n=3$ and $M$ is orientable, it follows from 
Proposition \ref{ynet} below that every elementary subgroup of $\Gamma$ is abelian. Hence this definition of Margulis number does generalize the one given above.

The {\it Margulis Lemma} \cite[Chapter D]{bp} implies that  there is a positive constant which is a Margulis number for
every hyperbolic $n$-manifold. The largest such number,
$\mu(n)$, is called the {\it Margulis constant} for hyperbolic
$n$-manifolds.

Meyerhoff showed in \cite{meyerhoff} that $\mu(3)>0.104$. Marc Culler has informed me that according to strong numerical evidence, 
$0.616$ fails to be a Margulis number for the hyperbolic $3$-manifold {\tt m027(-4,1)}, and hence 
$\mu(3) < 0.616$.

In \cite[Chapter D]{bp} it is explained in detail how a Margulis number $\mu$ for a hyperbolic manifold $M$ determines a canonical decomposition of $M$  into a {\it $\mu$-thin part}, consisting of cusp neighborhoods and tubes around closed geodesics,  and a {\it $\mu$-thick part} consisting of points where the injectivity radius is at least $\mu/2$. (In \cite{bp} the number $\mu$ is taken to be a Margulis {\it constant}, but the same arguments apply if $\mu$ is any Margulis number for $M$.

There are only finitely many topological possibilities for the $\mu$-thick part of $M$ given an upper bound on the volume of $M$. In the case $n=3$, the $\mu$-thick part is a $3$-manifold with torus boundary components,  and $M$ is diffeomorphic to the interior of a manifold obtained by Dehn fillings from its $\mu$-thick part.  

This makes estimation of the maximal Margulis number for $M$ a major step in understanding the geometric structure of $M$.  The larger $\mu$ is, the fewer possibilities there are for the $\mu$-thick part.

The proof of Theorem \ref{da goods} uses Proposition 8.5 of \cite{hakenmarg}, which in turn depends in an essential way on Theorem VI.4.1 of the Memoir \cite{js} by Bus Jaco and myself. This is an interesting example of how three-manifold topology can be applied to hyperbolic geometry. To deduce Theorem \ref{da goods} from \cite[Proposition 8.5]{hakenmarg} one needs results on algebraic convergence proved by T. Jorgensen, P. Klein and A. Weil. Many of the relevant definitions and results are reviewed in Section \ref{prelim} below.

I am grateful to Dick Canary for crucial help in arriving at the algebraic convergence arguments needed for the proof of Theorem \ref{da goods}; to Dave Futer for helping me to clarify the logical relationship between  Theorem \ref{da goods} and the results of \cite{hakenmarg}; to Al Marden for helping me locate the reference \cite{jk}; and to Marc Culler for computing the example discussed above.

\section{Preliminaries}\label{prelim}

As I mentioned in the introduction, a subgroup of $\isomplus(\HH^3)$ is said to be
{\it elementary} if it has an abelian subgroup of finite index.

In the case where $\Gamma$ is discrete and
torsion-free, I will say
 that the orientable hyperbolic $3$-manifold
$M=\HH^3/\Gamma$ is {\it elementary} if $\Gamma$ is elementary. If $M$ has finite volume then $M$ must be non-elementary. 

The following easy result gives a more direct characterization of elementary orientable hyperbolic $3$-manifolds.

\Proposition\label{ynet}
Every torsion-free, elementary, discrete subgroup of $\isomplus(\HH^3)$ is abelian.
\EndProposition

\Proof
Suppose that $\Gamma\le\isomplus(\HH^3)$ is  torsion-free, elementary, and discrete. For $1\ne\gamma\in\Gamma$, the subset $\Fix(\gamma) $ of the sphere at infinity is non-empty because $\langle\gamma\rangle$ is discrete and $\gamma$ has infinite order, and consists of at most two points because $\gamma\ne1$. For distinct elements $\gamma,\gamma'\in\Gamma$, we have $\Fix(\gamma)=\Fix(\gamma')$ if and only if $\gamma$ and $\gamma'$ commute. Now if $x$ and $y$ are non-trivial elements of $\Gamma$, the definition of an elementary group implies that $x^m$ and $y^n$ commute for some positive integers $m$ and $n$. We have
$$\Fix(x)=\Fix(x^m)=\Fix(y^n)=\Fix(y),$$
so that $x$ and $y$ commute.
\EndProof

We will need the following well-known fact:

\Proposition\label{rational polka dots}
For any compact, connected, orientable $3$-manifold  $N$, the total genus of $\partial N$, i.e. the sum of the genera of its
components, is bounded above by the first Betti number of $N$.
\EndProposition

\Proof
If one replaces the first Betti number of $N$, i.e. the dimension of the $\QQ$-vector space $H_1(N,\QQ)$, by the dimension of the $\ZZ_2$-vector space $H_1(N,\ZZ_2)$, the statement become Lemma 7.3 of \cite{last}. The proof of the latter result goes through without change if the coefficient field $\ZZ_2$ is replaced by $\QQ$.
\EndProof

The next result includes the second assertion of Theorem \ref{da goods}. The main ingredient in the proof is Proposition 8.5 of \cite{hakenmarg}, which asserts that $0.292$ is a Margulis number for every orientable $3$-manifold $M$ with $H_1(M;\QQ)\ne0$.

\Proposition\label{steven witt}
Every non-compact, orientable hyperbolic $3$-manifold admits $0.292$ as a Margulis number.
\EndProposition

\Proof
Let $M$ be a non-compact, orientable hyperbolic $3$-manifold.
write $M=\HH^3/\Gamma$ for some torsion-free
discrete subgroup $\Gamma$ of $\isomplus(\HH^3)$.
We must show that if $x$ and $y$ are non-commuting elements of $\Gamma$ 
and $P$ is a point in $\HH^3$, we have
$$\max(d(P,x\cdot P),d(P,y\cdot P))\ge0.29.$$

Set $\tGamma:=\langle x,y\rangle$. Then $\tM:=\HH^3/\tGamma$ is a hyperbolic manifold which covers $M$, and is therefore non-compact.
Since $\pi_1(\tM)\cong\tGamma=\langle x,y\rangle$ is finitely generated, $\tM$ has a {\it compact core} $N$ according to \cite{core}. By definition, $N\subset \tM$ is a compact submanifold, and the inclusion homomorphism $\pi_1(N)\to\pi_1(\tM)$ is an isomorphism.

Suppose that $\partial N$ has a sphere component $S$. Since $\tM$ is irreducible, $S$ bounds a ball $B\subset \tM$. If $B\supset N$ then $\tM$ is simply connected, a contradiction since $\tGamma \cong\pi_1(\tM)$ contains the non-commuting elements $x$ and $y$. Hence $B\cap N=S$, and it follows that $B\cap N$ is a compact core for $\tGamma$ with fewer boundary components than $N$. Hence if we choose $N$ among all compact cores so as to minimize its number of boundary components, then
$\partial N$ has no sphere component.

On the other hand, since $\tM$ is not closed, we have $\partial N\ne\emptyset$. Hence the total genus of $\partial N$, i.e. the sum of the genera of its components, is strictly positive. 

By Proposition \ref{rational polka dots}, the first Betti number of $N$ is bounded below by the total genus of $\partial N$. Hence $H_1(\tM;\QQ)\ne0$. According to \cite[Proposition 8.5]{hakenmarg}, this implies that $0.292$ is a Margulis number for
$\tM$. Since
$x$ and $y$ do not commute, we therefore have
$$\max(d(P,x\cdot P),d(P,y\cdot P))\ge0.292.$$
\EndProof

If $\Gamma$ is a group, we denote by $\isomplus(\HH^3)^{(\Gamma)}$ the set of all set-theoretical mappings of $\Gamma$ into $\isomplus(\HH^3)$. Thus $\isomplus(\HH^3)^{(\Gamma)}$ is the same as the product $\prod_{\gamma\in\Gamma}\isomplus(\HH^3)$, where $G_\gamma=\isomplus(\HH^3)$ for $\gamma\in\Gamma$. We shall always understand $\isomplus(\HH^3)$ to have its usual Lie group topology, and $\isomplus(\HH^3)^{(\Gamma)}$ to have the product topology. The set $\Hom(\Gamma,\isomplus(\HH^3))$ of all representations of $\Gamma$ in $\isomplus(\HH^3)$ is the subset of mappings $\rho\in \isomplus(\HH^3)^{(\Gamma)}$ which satisfy $\rho(\gamma\gamma')=\rho(\gamma)\rho(\gamma')$ for all $\gamma,\gamma'\in \isomplus(\HH^3)$. Hence $\Hom(\Gamma,\isomplus(\HH^3))$ is a closed subset of $\isomplus(\HH^3)^{(\Gamma)}$. It will be given the subspace topology. 

A sequence of representations $(\rho_j)_{j\in\NN}$ in $\Hom(\Gamma,\isomplus(\HH^3))$ is said to {\it converge algebraically} to a representation $\rho_\infty\in\Hom(\Gamma,\isomplus(\HH^3))$ if it converges to $\rho_\infty$ in the topology of $\Hom(\Gamma,\isomplus(\HH^3))$ described above.

I will denote by $D(\Gamma)$ the subset of $\Hom(\Gamma,\isomplus(\HH^3))$ consisting of all representations whose images are discrete and non-elementary, and by $DFC(\Gamma) \subset D(\Gamma)$ the subset of $\Hom(\Gamma,\isomplus(\HH^3))$ consisting of all faithful representations whose images are discrete and cocompact.

The following result summarizes the facts about algebraic convergence that will be needed in the next section.

\Theorem\label{jkw} Let $\Gamma$ be any group.
\begin{enumerate}
\item
 The set $D(\Gamma)$ is closed in $\Hom(\Gamma,\isomplus(\HH^3))$.
\item
If a sequence $(\phi_j)_{j\in\NN}$ in $D(\Gamma)$ converges
algebraically to a representation $\phi_\infty$, then for every sufficiently large $j$ there is a homomorphism $\psi_j:\phi_\infty(\Gamma)\to \phi_j(\Gamma)$ such that $\psi_j\circ\phi_\infty=\phi_j$.
\item
The set $DFC(\Gamma)$ is open in $\Hom(\Gamma,\isomplus(\HH^3))$ $DF(\Gamma)$.
\end{enumerate}
\EndTheorem

\Proof

In proving the first assertion we may assume that $D(\Gamma)\ne\emptyset$, so that $\Gamma$ is isomorphic to a discrete subgroup of $\isomplus(\HH^3)$ and is therefore countable. This implies that $\isomplus(\HH^3)^{(\Gamma)}$ is a first-countable space, and hence that $\Hom(\Gamma,\isomplus(\HH^3))$ is also first-countable. It therefore suffices to prove that every algebraic limit of representations in  $D(\Gamma)$ belongs to $D(\Gamma)$. But this is the first assertion of \cite[p. 326, Theorem]{jk}. 

The second assertion of Theorem \ref{jkw} is Theorem 2  of \cite{jor} (and is also the second assertion of \cite[p. 326, Theorem]{jk}).

The third assertion is the case $G=\isomplus(\HH^3)$ of the main theorem of \cite{weil}.
\EndProof

\section{Proof of the main theorem}\label{proof section}

The proof of Theorem \ref{da goods} will occupy this section. 

The second assertion of the theorem is immediate from Proposition \ref{steven witt}.

To prove the first assertion, we argue by contradiction.
Suppose that there is an infinite sequence $(M_j)_{j\ge0}$ of
pairwise non-isometric orientable hyperbolic $3$-manifolds
none of which admits $0.29$ as a Margulis number. For
each $j$ write $M_j=\HH^3/\Gamma_j$ for some torsion-free
discrete subgroup $\Gamma_j$ of $\isomplus(\HH^3)$.

For
each $j$, since $0.29$ is not a Margulis number for $M_j$,
there exist non-commuting elements
$x_j,y_j\in\Gamma_j$   and a point $P_j\in\HH^3$   such that 
$$\max(d(P_j,x_j\cdot P_j),d(P_j,y_j\cdot P_j))<0.29.$$
After replacing each $\Gamma_j$ by a suitable conjugate of itself in
$\isomplus(\HH^3)$, we may assume that the $P_j$ are all the same point of $\HH^3$,    which I will denote by $P$.    Thus for each $j$ we have
\begin{equation}\label{maximillian}\max(d(P,x_j\cdot P),d(P,y_j\cdot P))<0.29.\end{equation}

For each $j\in\NN$, the group $\tGamma_j:=\langle x_j,y_j\rangle$ is a subgroup of $\Gamma_j$, and is therefore discrete and torsion-free. Since $x_j$ and $y_j$ do not commute, it follows from Proposition \ref{ynet} that $\tGamma_j$ is non-elementary. Hence $\tM_j:=\HH^3/\tGamma_j$ is a non-elementary hyperbolic manifold which covers $M_j$.

It follows from (\ref{maximillian}) that the $x_j$ and $y_j$ all lie
in a compact subset of $\isomplus(\HH^3)$. Hence after passing to a
subsequence we may assume that $x_j\to x_\infty$ and $y_j\to y_\infty$
for some $x_\infty,y_\infty\in\isomplus(\HH^3)$. 

Let $F_2$ denote the free group on two generators $x$ and $y$. For
$1\le j\le\infty$, let us define a homomorphism
$\phi_j:F_2\to\isomplus(\HH^3)$ by $\phi_j(x)=x_j$ and
$\phi_j(y)=y_j$. The sequence $(\phi_j)_{j\in\NN}$ converges
algebraically to $\phi_\infty$. For each $j\in\NN$, since $\tGamma_j$ is discrete and non-elementary, we have $\phi_j\in D(\Gamma)$ in the notation of Section \ref{prelim}. Applying Assertion (1) of Theorem \ref{jkw}, with $\Gamma=F_2$ and $\rho_j=\phi_j$, we deduce that $\phi_j\in D(\Gamma)$; that is, $\tGamma_\infty:=\phi_\infty(F_2)$ is discrete and non-elementary. Furthermore, it follows from Assertion (2) of Theorem \ref{jkw} that for sufficiently large $j$ there is a homomorphism $\psi_j:\tGamma_\infty\to\tGamma_j$ such that $\psi_j\circ\phi_\infty=\phi_j$. After passing to a subsequence we may assume that such a $\psi_j$ exists for every $j\in\NN$. 

We regard $\psi_j$ as a representation of $\tGamma_\infty$ in $\isomplus(\HH^3)$, and we let $\psi_\infty:\tGamma_\infty\to \isomplus(\HH^3)$ denote the inclusion homomorphism. Since the sequence $(\phi_j)_{i\in\NN}$ converges algebraically to $\phi_\infty$, the
sequence $(\psi_j)_{i\in\NN}$ converges algebraically to $\psi_\infty$.

If $\gamma\in \tGamma_\infty$ has finite order, then for each $j\in\NN$ we have $\psi_j(\gamma)=1$, since $\tGamma_j$ is torsion-free. By algebraic convergence we have 
$$\gamma=\psi_\infty(\gamma)=\lim_{j\to\infty}\psi_j(\gamma)=1.$$
This shows that $\tGamma_\infty$ is
torsion-free. Hence $\tM_\infty=\HH^3/\tGamma_\infty$ acquires the structure of a hyperbolic $3$-manifold in a natural way.

We now distinguish two cases; in each case we shall obtain a contradiction.

{\bf Case I: $\tM_\infty$ is closed.} In this case $\Gamma_\infty$ is cocompact, so that in the notation of Section \ref{prelim} we have $\psi_\infty\in DFC(\Gamma_\infty)$. Since $(\psi_j)_{i\in\NN}$ converges algebraically to $\psi_\infty$, it then follows from Assertion (3) of Theorem \ref{jkw}, with $\Gamma=\tGamma_\infty$, that $\psi_j\in DFC(\tGamma_\infty)$ for every sufficiently large $j$. After passing to a subsequence we may therefore assume that $\psi_j\in DFC(\tGamma_\infty) $ for every $j\in\NN$. In particular $\psi_j:\tGamma_\infty\to\tGamma_j$  is an isomorphism for each $j$; and $\tM_j=\HH^3/\tGamma_j=\HH^3/\psi_j(\tGamma_\infty)$ is closed, and is therefore a finite-sheeted covering of $M_j$. For each $j\in\NN$, since $\tGamma_j$ is isomorphic to $\tGamma_\infty$, it follows from the Mostow rigidity theorem \cite[Chapter C]{bp}  that $\tM_j$ is isometric to $\tM_\infty$.

Let $\Delta $ denote the diameter of $ \tM_\infty$. Then since each $M_j$ has a
finite-sheeted covering isometric to $\tM_\infty$, each $M_j$ has diameter
at most $\Delta $.

Let $v$ denote the infimum of the volumes of all
closed hyperbolic $3$-manifolds;  we have $v>0$, for example by
\cite[Theorem 1]{meyerhoff}. Since each $M_j$ has volume at least $v$,
diameter at most $\Delta$, constant curvature $-1$ and dimension $3$, it
follows from the main theorem of \cite{peters} that the $M_j$
represent only finitely many diffeomorphism types. By the Mostow
rigidity theorem, they represent only finitely many isometry
types. This gives the required contradiction in Case I.

{\bf Case II: $\tM_\infty$ is not closed.} In this case it follows from Proposition \ref{steven witt} that $0.292$ is a Margulis number for
$\tM_\infty$. Since $\tGamma_\infty=\langle x_\infty,y_\infty\rangle$ is non-elementary, 
$x_\infty$ and $y_\infty$ do not commute. We therefore have
$$\max(d(P,x_\infty\cdot P),d(P,y_\infty\cdot P))\ge0.292.$$
On the other hand, it follows from
(\ref{maximillian}), upon taking limits as $j\to\infty$, that 
$$\max(d(P,x_\infty\cdot P),d(P,y_\infty\cdot P))\le0.29.$$
This gives the required contradiction in Case II.

\bibliographystyle{plain}
\bibliography{finiteness}

\end{document}